\documentstyle[amscd,amssymb,verbatim,12pt]{amsart}
\newcommand{\bb}{\mbox{\rm b}}

\newcommand{\C}{{\Bbb C}}

\newcommand{\HH}{\mbox{\rm H}}

\newcommand{\Id}{\mbox{\rm Id}}

\newcommand{\Image}{\mbox{\rm Im}}

\newcommand{\Ker}{\mbox{\rm Ker}}

\newcommand{\R}{{\Bbb R}}

\newcommand{\vol}{\mbox{\rm vol}}
\newcommand{\Z}{{\Bbb Z}}
\theoremstyle{plain}

\newtheorem{lemma}{Lemma}
\newtheorem{theorem}{Theorem}
\newtheorem{proposition}{Proposition}

\numberwithin{equation}{section}
\renewcommand{\rm}{\normalshape}
\begin{document}
\title{Remark About the Spectrum of the $p$-Form Laplacian Under a 
Collapse with Curvature Bounded Below}
\author{John Lott}
\address{Department of Mathematics\\
University of Michigan\\
Ann Arbor, MI  48109-1109\\
USA}
\email{lott@@math.lsa.umich.edu}
\thanks{Research supported by NSF grant DMS-0072154}
\subjclass{Primary: 58G25; Secondary: 53C23}
\date{February 11, 2002}
\maketitle
\begin{abstract}
We give a lower bound on the number of small positive eigenvalues of the
$p$-form Laplacian in a certain type of collapse with curvature bounded
below.
\end{abstract}

\section{Introduction} \label{Introduction}

A general problem in spectral geometry is to estimate the eigenvalues
of the $p$-form Laplacian on a closed Riemannian manifold $M$ in terms
of the geometry of $M$. 
From Hodge theory, the number of zero eigenvalues is $\bb_p(M)$, the
$p$-th Betti number of $M$. Hence the issue is to understand the
positive eigenvalues.  The papers
\cite{Colbois-Courtois (1990)}, \cite{Lott (2001)} and \cite{Lott (2002)}
study the case when one assumes an upper bound on the diameter of the
manifold and double-sided bounds on the sectional curvatures.
An important phenomenon
is the possible appearance of positive
eigenvalues of the $p$-form Laplacian that approach zero
as a manifold collapses with bounded curvature.

The analysis of 
\cite{Lott (2001)} and \cite{Lott (2002)}
uses the results of Cheeger, Fukaya and Gromov on the geometric structure of
manifolds that collapse with double-sided curvature bounds.  If one only 
assumes a lower sectional curvature bound then there are some structure
results about collapsing in \cite{Fukaya-Yamaguchi (1992)} and
\cite{Yamaguchi (1991)},
but the theory is less developed than in the bounded curvature case.

In this paper we look at the small positive eigenvalues of the
$p$-form Laplacian in an example of collapse with
curvature bounded below.
Namely, suppose that a compact
Lie group $G$ acts isometrically on $M$ on the left. Give $G$ a left-invariant
Riemannian metric. For $\epsilon \: > \: 0$, let $\epsilon G$ denote $G$
with its Riemannian metric multiplied by $\epsilon^2$. Let $M_\epsilon$
denote $M \: = \: G \backslash (\epsilon G \times M)$ 
equipped with the quotient Riemannian metric $g_\epsilon$,
where $G$ acts
diagonally on $\epsilon G \times M$ on the left. If $G$ is connected then
$\lim_{\epsilon \rightarrow 0} M_\epsilon \: = \: G \backslash M$ 
in the Gromov-Hausdorff topology,
and as $\epsilon$ goes to zero, the sectional curvatures of $M_\epsilon$ 
stay uniformly bounded below \cite{Yamaguchi (1991)}.

For notation,
if $M$ is a smooth connected closed manifold with Riemannian metric $g$,
let $\{ \lambda_{p,j}(M, g) \}_{j=1}^\infty$ denote the
eigenvalues (counted with multiplicity) of the Laplacian on
$\overline{\Image(d)} \subset \Omega^p_{L^2}(M)$.

\begin{theorem} \label{Theorem 1} 
If $j \: = \:
\dim \left( \Ker ( \HH^p(G \backslash M; \R) \rightarrow \HH^p(M; \R)) \right)$
then $\lim_{\epsilon \rightarrow 0} \lambda_{p,j}(M_\epsilon,
g_\epsilon) \: = \: 0$.
\end{theorem} 

In Section \ref{Proof} we prove Theorem \ref{Theorem 1}.  
The main points of the proof
are the use of a certain variational expression for $\lambda_{p,j}(M, g)$, 
due to Cheeger and Dodziuk \cite{Dodziuk (1982)}, and the avoidance of
dealing with the detailed orbit structure of the group action.
We then look at the example of an $S^1$-action on $S^{2n}$,
which is the suspension of the Hopf action of $S^1$ on $S^{2n-1}$,
and show that our results slightly improve those of
Takahashi \cite{Takahashi (2001)}.  
In Section \ref{Remarks} we make some further remarks.

I thank Junya Takahashi for sending me a copy of his paper.

\section{Proof of Theorem \ref{Theorem 1}} \label{Proof}

Let ${\frak g}$ be the Lie algebra of $G$. It acquires an inner product from
the left-invariant Riemannian metric on $G$.  Given ${\frak x} \in {\frak g}$,
let ${\frak X}$ be the corresponding vector field on $M$. 
Let $i_{\frak X}$ denote interior multiplication by ${\frak X}$.

Let $\Omega^*(M)$ denote the smooth differential forms on $M$.
Let $\Omega^*_{L^2}(M)$ be the $L^2$-completion of $\Omega^*(M)$.
Put
\begin{equation} \label{2.1}
\Omega^*_{max}(M) \: = \: \{ \omega \in \Omega^*_{L^2}(M) \: : \:
d\omega \in \Omega^{*+1}_{L^2}(M) \},
\end{equation}
where $d\omega$ is originally defined distributionally.

Put
\begin{equation} \label{2.2}
\Omega^*_G(M) \: = \: \{ \omega \in \Omega^*(M) \: : \: g \cdot \omega \: = \:
\omega \text{ for all } g \in G \} 
\end{equation}
and
\begin{equation} \label{2.3}
\Omega^*_{basic}(M) 
\: = \: \{ \omega \in \Omega^*_G(M) \: : \: i_{\frak X} \omega \: = \: 0
\text{ for all } {\frak x} \in {\frak g} \}.
\end{equation}
Let $\Omega^*_{G, L^2}(M)$ and 
$\Omega^*_{basic, L^2}(M)$ be the $L^2$-completions of 
$\Omega^*_{G}(M)$ and 
$\Omega^*_{basic}(M)$, respectively.
Put
\begin{equation} \label{2.4}
\Omega^*_{basic,max}(M) \: = \: \{ \omega \in \Omega^*_{basic, L^2}(M) \: : \:
d\omega \in \Omega^{*+1}_{basic, L^2}(M) \},
\end{equation}
where $d\omega$ is originally defined distributionally.
Then $\Omega^*_{basic,max}(M)$ is a complex.

From \cite{Koszul (1953)} and \cite{Verona (1988)}, 
the cohomology of the complex $\Omega^*_{basic}(M)$
is isomorphic to $\HH^*(G \backslash M; \R)$.
\begin{lemma} \label{Lemma 1}
The cohomology of the complex $\Omega^*_{basic,max}(M)$
is isomorphic to $\HH^*(G \backslash M; \R)$.
\end{lemma}
\begin{pf}
The proof is essentially the same as that of \cite{Verona (1988)}.
The only point to note is that the homotopy operator $A$ used
in the Poincar\'e lemma in \cite{Verona (1988)} 
sends $\Omega^*_{basic,max}$ to itself.
\end{pf}

The quotient map $p \: : \: \epsilon G \: \times \: M \rightarrow
M_\epsilon$ defines a principal $G$-bundle. Pullback gives an isomorphism
$p^* \: : \:
\Omega^*(M_\epsilon) \: \cong \: \Omega^*_{basic}(\epsilon G \: \times \: M)$.
The parallelism of $G$ gives an isomorphism
\begin{equation} \label{2.5}
\Omega^*(\epsilon G \: \times \: M) \: \cong \:
C^\infty(G) \: \otimes \: \Lambda^*({\frak g}^*) \: \otimes \: \Omega^*(M).
\end{equation}
Taking $G$-invariants gives isomorphisms
\begin{equation} \label{2.6}
\Omega^*_G(\epsilon G \: \times \: M) \: \rightarrow \:
\left(
C^\infty(G) \: \otimes \: \Lambda^*({\frak g}^*) \: \otimes \: \Omega^*(M)
\right)^G \: \stackrel{\beta}{\rightarrow} \: \Lambda^*({\frak g}^*) \:
\otimes \: \Omega^*(M),
\end{equation}
where $\beta$ comes from the map which
sends $\sum_k f_k \: \otimes \: \eta_k \: 
\otimes \: \omega_k
\in  
C^\infty(G) \: \otimes \: \Lambda^*({\frak g}^*) \: \otimes \: \Omega^*(M)$ to
$\sum_i f_k(e) \: \eta_k \: \otimes \: \omega_k \in 
\Lambda^*({\frak g}^*) \: \otimes \: \Omega^*(M)$.

Let $\{ {\frak x}_j \}_{j=1}^{dim(G)}$ be a basis of ${\frak g}$.
For ${\frak x} \in {\frak g}$, let $e({\frak x}^*)$ denote exterior
multiplication by ${\frak x}^*$ on $\Lambda^*({\frak g}^*)$.

\begin{lemma} \label{Lemma 2}
There is an isomorphism of complexes ${\cal I} \: : \: \Omega^*(M) \rightarrow
\Omega^*_{basic}(\epsilon G \: \times \: M) \subset
\Lambda^*({\frak g}^*) \: \otimes \: \Omega^*(M)$ given by 
\begin{align} \label{2.7}
{\cal I}(\sigma) \: & = \:
\left( \prod_{j=1}^{dim(G)} \left( 1 \: - \: e({\frak x}_j^*) \: \otimes \: 
i_{{\frak X}_j} \right) \right) \: (1 \: \otimes \: \sigma) \notag \\
& = \:
\sum_{k=0}^{dim(G)} (-1)^k \: 
\sum_{1 \: \le \: j_1 \: < \: \ldots \: < \:
j_k \le
dim(G)} ({\frak x}_{j_k}^* \wedge \ldots \wedge {\frak x}_{j_1}^*) \: \otimes
\: i_{{\frak X}_{j_1}} \: \ldots i_{{\frak X}_{j_k}} \sigma.
\end{align}
\end{lemma}
\begin{pf}
If $\sum_k f_k \: \otimes \: \eta_k \: \otimes \: \omega_k \in
\left(
C^\infty(G) \: \otimes \: \Lambda^*({\frak g}^*) \: \otimes \: \Omega^*(M)
\right)^G$ is $G$-basic
then for ${\frak x} \in {\frak g}$, we also have
\begin{equation} \label{2.8}
\sum_k \left( f_k 
\: \otimes \: i_{\frak x} \eta_k \: \otimes \: \omega_k \: + \:
(-1)^{|\eta_k|} \:  
f_k \: \otimes \: \eta_k \: \otimes \: i_{\frak X} \omega_k \right) 
\: = \: 0.
\end{equation} 
Then
\begin{equation} \label{2.9}
\sum_k \left( f_k(e) 
\: i_{\frak x} \eta_k \: \otimes \: \omega_k \: + \: (-1)^{|\eta_k|} \:  
f_k(e) \: \eta_k \: \otimes \: i_{\frak X} \omega_k \right) 
\: = \: 0,
\end{equation}
i.e. if $\sum_k \eta_k \: \otimes \: \omega_k$ lies in the image of
$\beta$ restricted to $\Omega^*_{basic}(\epsilon G \: \times \: M)$ then
\begin{equation} \label{2.10}
\sum_k \left( i_{\frak x} \eta_k \: \otimes \: \omega_k \: + \: 
(-1)^{|\eta_k|} \: 
\eta_k \: \otimes \: i_{\frak X} \omega_k \right) 
\: = \: 0.
\end{equation}
It follows that $\sum_k \eta_k \: \otimes \: \omega_k$ can be written as
${\cal I}(\sigma)$ for some $\sigma \in \Omega^*(M)$. Thus ${\cal I}$ is
surjective.  It is clearly injective.

It remains to show that ${\cal I}$ is a morphism of complexes.  Let
$d^{inv}$ denote the (finite-dimensional) differential on
$\Lambda^*({\frak g}^*)$.
If an element of $\Omega^*_{G}(\epsilon G \: \times \: M)$ is 
represented as $\sum_k f_k \: \otimes \: \eta_k \: 
\otimes \: \omega_k
\in  
C^\infty(G) \: \otimes \: \Lambda^*({\frak g}^*) \: \otimes \: \Omega^*(M)$
then the $G$-invariance implies that
for ${\frak x} \in {\frak g}$,
\begin{equation} \label{2.11}
\sum_k \left( {\frak x} f_k 
\: \otimes \: \eta_k \: \otimes \: \omega_k \: + \: 
f_k \: \otimes \: \eta_k \: \otimes \: {\cal L}_{\frak X} \omega_k \right) 
\: = \: 0.
\end{equation}  

The differential of $\sum_k f_k \: \otimes \: \eta_k \: 
\otimes \: \omega_k$ is represented by  
\begin{equation} \label{2.12}
\sum_k \left( \sum_{j=1}^{dim(G)} 
{\frak x}_j f_k \: \otimes \: e({\frak x}_j^*) \eta_k \: 
\otimes \: \omega_k \: + \:
f_k \: \otimes \: d^{inv} \eta_k \: 
\otimes \: \omega_k \: + \: (-1)^{|\eta_k|} \:
f_k \: \otimes \: \eta_k \: 
\otimes \: d\omega_k \right).
\end{equation}
From (\ref{2.11}), this equals
\begin{equation} \label{2.13}
\sum_k \left( - \: \sum_{j=1}^{dim(G)} 
f_k \: \otimes \: e({\frak x}_j^*) \eta_k \: 
\otimes \: {\cal L}_{{\frak X}_j} \omega_k \: + \:
f_k \: \otimes \: d^{inv} \eta_k \: 
\otimes \: \omega_k \: + \: (-1)^{|\eta_k|} \:
f_k \: \otimes \: \eta_k \: 
\otimes \: d\omega_k \right).
\end{equation}
Using $\beta$, it follows that the induced differential on 
$\Lambda^*({\frak g}^*) \:
\otimes \: \Omega^*(M)$ sends 
$\sum_k \eta_k \: 
\otimes \: \omega_k$ to
\begin{equation} \label{2.14}
\sum_k \left( - \: \sum_{j=1}^{dim(G)} e({\frak x}_j^*) \eta_k \: 
\otimes \: {\cal L}_{{\frak X}_j} \omega_k \: + \:
d^{inv} \eta_k \: 
\otimes \: \omega_k \: + \: (-1)^{|\eta_k|} \: 
\eta_k \: 
\otimes \: d\omega_k \right).
\end{equation}
One can check that when this acts on ${\cal I}(\sigma)$, the result is
${\cal I}(d\sigma)$. Thus ${\cal I}$ is an isomorphism of complexes.
\end{pf}

In fact, under our identifications, ${\cal I}$ is the same as $p^*$.

Let $M^{reg}$ be the union of the
principal orbits for the $G$-action on $M$. It is a dense open subset
of $M$ with full measure. If $m \in M^{reg}$, 
let $H \subset G$ 
be its isotropy subgroup, with Lie algebra ${\frak h}$.
Define ${\alpha} 
\: : \: {\frak g} \rightarrow
T_mM$ by ${\alpha}({\frak x}) 
\: = \: {\frak X}_m$. It passes to an injection
$\overline{\alpha} 
\: : \: {\frak g}/{\frak h} \rightarrow
T_mM$.  For $\epsilon \: \ge 0$, put $\rho_\epsilon(m) \: = \: 
\det^{1/2}(\epsilon^2 \: \Id. \big|_{{\frak g}/{\frak h}} 
\: + \: \overline{\alpha}^* \: \overline{\alpha})$.
If $m \notin M^{reg}$,
put $\rho_\epsilon(m) \: = \: 0$. Note that for $\epsilon \: > \: 0$,
$\rho_\epsilon^{-1}(m) \: < \: \rho^{-1}_0(m)$.
\begin{lemma} \label{Lemma 3}
$\rho^{-1}_0 \in L^1(M, d\vol)$.
\end{lemma}
\begin{pf}
If $m \in M^{reg}$ then
up to an overall constant, $\rho_0(m)$ is the volume of the orbit
$G \cdot m$. Then $\int_{M^{reg}} \rho_0^{-1}(m) \: d\vol(m)$ is 
proportionate to the volume of
$G \backslash M^{reg} \subset G \backslash M$, which is seen to be finite.
\end{pf}

Let $\{ {\frak x}_j \}_{j=1}^{dim(G)}$ be an orthonormal basis of ${\frak g}$.

\begin{lemma} \label{Lemma 4}
For $\epsilon \: > \: 0$, there is a positive constant $C(\epsilon)$ such that
$\Omega^*(M_\epsilon)$ 
is isometrically isomorphic to $\Omega^*(M)$ with the new norm
\begin{align} \label{2.15}
\parallel \omega \parallel_\epsilon^2 \: = \: & C(\epsilon) \: 
\int_M \rho_\epsilon^{-1}(m) \notag \\
& \left( | \omega(m) |_M^2 \: + \:
\sum_{k \: = \: 1}^{dim(G)} \: \epsilon^{- \: 2k} \:
\sum_{1 \: \le \: j_1 \: < \: \ldots \: < \:
j_k \le
dim(G)} | i_{{\frak X}_{j_1}} \: \ldots i_{{\frak X}_{j_k}} \omega(m)
|_M^2 \right) \: d\vol(m).
\end{align}
\end{lemma}
\begin{pf}
We can compute the norm squared of 
$\omega \in \Omega^*(M_\epsilon)$ by taking the local norm squared of
$p^* \omega$ on $\epsilon G \: \times M^{reg}$, dividing by the function
which assigns to $(g,m) \in \epsilon G \: \times M^{reg}$ the volume
of the orbit $G \cdot (g, m)$, 
and integrating over $\epsilon G \: \times M^{reg}$. If $m \in M^{reg}$ then
the relative volume of $G \cdot (g, m)$ is 
\begin{equation} \label{2.16}
{\det}^{1/2} (\epsilon^2 \: \Id. \big|_{{\frak g}} 
\: + \: {\alpha}^* {\alpha}) \: = \:
\epsilon^{dim(H)} \: 
\rho_\epsilon(m).
\end{equation} 
The map $\beta$ of (\ref{2.6}) is an
isometry, up to a constant. As $\{ \epsilon^{-1} \: 
{\frak x}_j \}_{j=1}^{dim(G)}$ is
an orthonormal basis for $T_e(\epsilon G)$, the lemma follows from
Lemma \ref{Lemma 2}.
\end{pf}
{\bf Proof of Theorem \ref{Theorem 1} : } 
Put $\lambda_{p,j}(\epsilon) \: = \: \lambda_{p,j}(M_\epsilon, g_\epsilon)$.
From \cite{Dodziuk (1982)},
\begin{equation} \label{2.17}
\lambda_{p,j}(\epsilon) \: = \: \inf_V \: \sup_{\eta \in V - 0} \:
\sup_{\theta \in d^{-1}(\eta)} \:
\frac{\parallel \eta \parallel_\epsilon^2}{\parallel 
\theta \parallel_\epsilon^2},
\end{equation}
where $V$ ranges over $j$-dimensional subspaces of $\Image \left( d \: : 
\Omega^{p-1}(M) \rightarrow 
\Omega^p(M) \right)$, and $\theta \in d^{-1}(\eta) \subset
\Omega^{p-1}(M)$.

Take $j \: = \:
\dim \left( \Ker ( \HH^p(G \backslash M; 
\R) \rightarrow \HH^p(M; \R)) \right)$.
From Lemma \ref{Lemma 1}, the inclusion of complexes
$\Omega^*_{basic}(M) \rightarrow \Omega^*_{basic, max}(M)$ induces an
isomorphism on cohomology. Then there is 
a $j$-dimensional subspace $V$ of 
\begin{equation} \label{2.18}
\Ker \left( d \: : \: \Omega^p_{basic}(M) \rightarrow \Omega^{p+1}_{basic}(M)
\right) \: \cap \: 
\Image \left( d \: : \: \Omega^{p-1}(M) \rightarrow \Omega^{p}(M)
\right)
\end{equation}
such that if $\eta \in V - 0$ then 
$\eta \notin \Image
\left( d \: : \: \Omega^{p-1}_{basic,max}(M) \rightarrow 
\Omega^{p}_{basic, L^2}(M)
\right)$. We claim that 
\begin{equation} \label{2.19}
\lim_{\epsilon \rightarrow 0} \: 
\sup_{\eta \in V - 0} \:
\sup_{\theta \in d^{-1}(\eta)} \:
\frac{\parallel \eta \parallel_\epsilon^2}{\parallel 
\theta \parallel_\epsilon^2} \: = \: 0.
\end{equation}
This will suffice to prove the theorem.

Suppose that (\ref{2.19}) is not true.  Then there are a constant
$c \: > \: 0$, a sequence
$\{\epsilon_r \}_{r=1}^\infty$ in $\R^+$ approaching zero, a sequence
$\{\eta_r\}_{r=1}^\infty$ in $V - 0$ and a sequence
$\{\theta_r\}_{r=1}^\infty$ in $\Omega^{p-1}(M)$ such that
for all $r$, $d \theta_r \: = \: \eta_r$ and
\begin{equation} \label{2.20}
\frac{\parallel \eta_r \parallel_{\epsilon_r}^2}{\parallel 
\theta_r \parallel_{\epsilon_r}^2} \: \ge \: c.
\end{equation}
Doing a Fourier decomposition of $\theta_r$ with respect to $G$,
the ratio in (\ref{2.20}) does not decrease
if we replace $\theta_r$ by its $G$-invariant component.
Thus we may assume that $\theta_r$ is $G$-invariant.

Without loss of generality, we can replace the norm $\parallel \cdot 
\parallel_\epsilon$ of (\ref{2.15}) by the same norm divided by $C(\epsilon)$,
which we again denote by $\parallel \cdot 
\parallel_\epsilon$.
As $\eta_r$ is smooth on $M$, it follows from Lemma \ref{Lemma 3} that
the function $\rho_0^{-1}(m) \: | \eta_r(m) |_M^2$ is integrable
on $M$. Without loss of generality, we may assume that
\begin{equation} \label{2.21}
\int_M \rho_0^{-1}(m) \: | \eta_r(m) |_M^2 \: d\vol(m) \: = \:
1.
\end{equation}
As $\{\eta_r\}_{r=1}^\infty$ lies in the sphere of a finite-dimensional space,
there will be a subsequence, which we relabel as $\{\eta_r\}_{r=1}^\infty$,
that converges smoothly to some $\eta_\infty \in V - 0$. 

From (\ref{2.20}),
\begin{align} \label{2.22}
\parallel \theta_r \parallel_{\epsilon_r}^2 \: & \le \: c^{-1}
\parallel \eta_r \parallel_{\epsilon_r}^2 \: = \:
c^{-1} \: \int_M \rho_{\epsilon_r}^{-1}(m) \: | \eta_r(m) |_M^2 \: 
d\vol(m) \notag \\ 
& \le \: c^{-1} \: 
\int_M \rho_0^{-1}(m) \: | \eta_r(m) |_M^2 \: d\vol(m) \: = \:
c^{-1}.
\end{align}
For large $r$,
\begin{equation} \label{2.23}
\int_M | \theta_r(m) |_M^2 \: d\vol(m) \: \le \: 
(\inf_M \rho_{\epsilon_r}^{-1})^{-1} \: \int_M \rho_{\epsilon_r}^{-1}(m) \: | 
\theta_r(m) |_M^2 \: d\vol(m) \: \le \: (\inf_M \rho_1^{-1})^{-1} \: c^{-1}.
\end{equation}

We now work with respect to the metric $g$ on $M$.
By weak-compactness of the unit ball in $L^2$, 
there is a subsequence of $\{ \theta_r \}_{r=1}^\infty$, which we
relabel as $\{ \theta_r \}_{r=1}^\infty$, that converges weakly in $L^2$ to
some $\theta_\infty \in \Omega^{p-1}_{G, L^2}(M)$. 
Then for $\sigma \in \Omega^p(M)$,
\begin{equation} \label{2.24}
\langle \sigma, \eta_\infty 
\rangle_M \: - \: \langle d^* \sigma, \theta_\infty \rangle_M \: =
\: \lim_{r \rightarrow \infty} \left( 
\langle \sigma, \eta_r
\rangle_M \: - \: \langle d^* \sigma, \theta_r \rangle_M \right) \: =
\: \lim_{r \rightarrow \infty} 
\langle \sigma, \eta_r \: - \: d \theta_r \rangle_M \: = \: 0.
\end{equation}
Thus $\theta_\infty \in \Omega^{p-1}_{max}(M)$ and $d\theta_\infty \: = \:
\eta_\infty$. 

From (\ref{2.22}), 
we also obtain that for each $1 \: \le \: j \: \le \: \dim(G)$,
\begin{equation} \label{2.25}
\int_M | i_{{\frak X}_j} \theta_r(m) |_M^2 \: d\vol(m) \: \le \: 
(\inf_M \rho_{\epsilon_r}^{-1})^{-1} \: \int_M \rho_{\epsilon_r}^{-1}(m) \: | 
i_{{\frak X}_j} \eta_r(m) |_M^2 \: d\vol(m) \: 
\le \: (\inf_M \rho_1^{-1})^{-1} \: c^{-1} \: \epsilon_r^2.
\end{equation}
Then for all $\sigma \in \Omega^{p-2}(M)$,
\begin{equation} \label{2.26}
\langle \sigma, i_{{\frak X}_j} \theta_\infty \rangle_M \: = \:
\langle (i_{{\frak X}_j})^* \sigma, \theta_\infty \rangle_M \: = \: 
\lim_{r \rightarrow \infty}
\langle (i_{{\frak X}_j})^* \sigma, \theta_r \rangle_M \: = \: 
\lim_{r \rightarrow \infty}
\langle  \sigma, i_{{\frak X}_j} \theta_r \rangle_M \: = \: 0. 
\end{equation}
Thus $i_{{\frak X}_j} \theta_\infty \: = \: 0$ and
$\theta_\infty \in \Omega^{p-1}_{basic,max}(M)$.
Hence $\eta_\infty \in \Image
\left( d \: : \: \Omega^{p-1}_{basic,max}(M) \rightarrow 
\Omega^{p}_{basic, L^2}(M)
\right)$, which is a contradiction. $\square$ \\ \\
{\bf Example : } Let $G \: = \: U(1)$ act on $M \: = \: S^{2n}$ by the
suspension of the Hopf action of $U(1)$ on $S^{2n-1}$. Then
$G \backslash M 
\: = \: U(1) \backslash S^{2n}$ is the suspension of $\C P^{n-1}$.
One finds that 
$\Ker ( \HH^p(G \backslash M; \R) 
\rightarrow \HH^p(M; \R))$ is nonzero if and only if
$p \: \in \: \{ 3, 5, \ldots, 2n-1 \}$. From Theorem \ref{Theorem 1}, 
as $\epsilon
\rightarrow 0$ there are small eigenvalues of the $p$-form Laplacian on
$\overline{\Image(d)} \subset \Omega^p_{L^2}(M_\epsilon)$ for
$p \: \in \: \{ 3, 5, \ldots, 2n-1 \}$.
From the Hodge decomposition, there will also be 
small eigenvalues of the $p$-form Laplacian on
$\overline{\Image(d^*)} \subset \Omega^p_{L^2}(M_\epsilon)$ for
$p \: \in \: \{ 2, 4, \ldots, 2n-2 \}$. Then using Hodge duality,
one concludes that there are small eigenvalues on \\
1. $\overline{\Image(d^*)} \subset \Omega^1_{L^2}(M_\epsilon)$, \\
2. $\overline{\Image(d)} \subset \Omega^p_{L^2}(M_\epsilon)$ and
$\overline{\Image(d^*)} \subset \Omega^p_{L^2}(M_\epsilon)$ for
$p \in \{ 2, 3, 4, \ldots, 2n-3, 2n-2 \}$, and \\
3. $\overline{\Image(d)} \subset \Omega^{2n-1}_{L^2}(M_\epsilon)$. \\
This slightly sharpens \cite[Theorem 1.2]{Takahashi (2001)}. Note that
from eigenvalue estimates for the scalar Laplacian \cite{Berard (1988)},
there are no small eigenvalues on 
$\overline{\Image(d^*)} \subset \Omega^0_{L^2}(M_\epsilon)$,
$\overline{\Image(d)} \subset \Omega^1_{L^2}(M_\epsilon)$,
$\overline{\Image(d^*)} \subset \Omega^{2n-1}_{L^2}(M_\epsilon)$ or
$\overline{\Image(d)} \subset \Omega^{2n}_{L^2}(M_\epsilon)$.

\section{Remarks} \label{Remarks}
\noindent
1. In the case of a locally-free torus action, there is some intersection
between Theorem \ref{Theorem 1} and the results of 
\cite{Colbois-Courtois (1990)},
\cite{Lott (2001)} and \cite{Lott (2002)}. 
In \cite{Lott (2002)} one deals with the cohomology of a certain
$\Z$-graded sheaf
$\HH^*(A^\prime_{[0]})$ on the limit space $X$. In the case of a collapsing
coming from a locally-free torus action,
Theorem \ref{Theorem 1} 
is a statement about the case $* \: = \: 0$, when the sheaf
$\HH^0(A^\prime_{[0]})$ is the constant $\R$-sheaf on $X$.
Of course, the result of
Theorem \ref{Theorem 1} will generally not give 
all of the small positive eigenvalues that arise in a collapse.
As seen in the Example, one can obtain more small
eigenvalues just from Hodge duality. \\ \\
2. Theorem \ref{Theorem 1} indicates that the relevant cohomology of the
limit space is the ordinary cohomology, as opposed for example to the
$L^2$-cohomology. This is consistent with the results of 
\cite{Lott (2002)} in the bounded curvature case.\\ \\
3. If $G$ has positive dimension and acts effectively on $M$ then
Theorem \ref{Theorem 1} describes small positive eigenvalues 
in a collapsing situation.  In some
noncollapsing situations, one can show that small
eigenvalues do not exist.  Here is one such criterion.

\begin{proposition} \label{Proposition 1}
 Let ${\cal M}$ be a collection of closed $n$-dimensional
Riemannian manifolds, with $n \: > \: 0$.
Give ${\cal M}$ the Lipschitz metric, coming from biLipschitz homeomorphisms.
Suppose that ${\cal M}$ can be covered by a finite number of metric balls.
For $p \in \Z \cap [0,n]$ and $j \: \ge \: 0$,
there
are positive numbers $a_{p,j}$ and $A_{p,j}$ so that if
$(M, g) \in {\cal M}$ then
$a_{p,j} \: \le \: \lambda_{p,j}(M, g) \: \le \: A_{p,j}$, and
$\lim_{j \rightarrow \infty} a_{p,j} \: = \: \infty$.
\end{proposition}
\begin{pf}
Suppose first
that for some $p$ and $j$, there is no upper bound on
$\lambda_{p,j}(M, g)$ as 
$(M, g)$ ranges over ${\cal M}$.
Then there is a 
sequence $\{(M_i, g_i)\}_{i=1}^\infty$ in ${\cal M}$ with the property
that $\lim_{i \rightarrow \infty} \lambda_{p,j}(M_i, g_i) \: = \: \infty$.
A subsequence of $\{(M_i, g_i)\}_{i=1}^\infty$, which we relabel as
$\{(M_i, g_i)\}_{i=1}^\infty$, will have finite distance from some
$(M_\infty, g_\infty) \in {\cal M}$. Then there are a number 
$\epsilon \: \ge \: 0$ and a  
sequence of biLipschitz homeomorphisms
$h_i \: : \: M_\infty \rightarrow M_i$ so that for all $i$,
\begin{equation} \label{3.1}
e^{- \: \epsilon} \: g_\infty \: \le \: h_i^* g_i \: \le \: 
e^{\epsilon} \: g_\infty.
\end{equation}
Here $h_i^* g_i$ is a Lipschitz metric on $M_\infty$.
From Hodge theory,
\begin{equation} \label{3.2}
\lambda_{p,j}(M_i, g_i) \: = \: \inf_V \: \sup_{\eta \in V - 0} \:
\sup_{\theta \in d^{-1}(\eta)} \:
\frac{\parallel \eta \parallel_{M_i}^2}{\parallel 
\theta \parallel_{M_i}^2},
\end{equation}
where $V$ ranges over $j$-dimensional subspaces of $\Image \left( d \: : 
\Omega^{p-1}_{max}(M_i) \rightarrow 
\Omega^p_{L^2}(M_i) \right)$, and $\theta \in d^{-1}(\eta) \subset
\Omega^{p-1}_{max}(M_i)$. By naturality,
\begin{equation} \label{3.3}
\lambda_{p,j}(M_i, g_i) \: = \: \inf_V \: \sup_{\eta \in V - 0} \:
\sup_{\theta \in d^{-1}(\eta)} \:
\frac{\parallel \eta \parallel_{h_i^* g_i}^2}{\parallel 
\theta \parallel_{h_i^* g_i}^2},
\end{equation}
where $V$ ranges over $j$-dimensional subspaces of $\Image \left( d \: : 
\Omega^{p-1}_{max}(M_\infty) \rightarrow 
\Omega^p_{L^2}(M_\infty) \right)$, and $\theta \in d^{-1}(\eta) \subset
\Omega^{p-1}_{max}(M_\infty)$.

As in \cite{Dodziuk (1982)}, it follows 
from (\ref{3.1}) and (\ref{3.3}) that there is a positive integer 
$J$ which only depends on $n$ so that 
\begin{equation} \label{3.4}
e^{- \: J \epsilon} \: \lambda_{p,j}(M_\infty, g_\infty) \: \le \: 
\lambda_{p,j}(M_i, g_i) \: \le \: 
e^{J \epsilon} \: \lambda_{p,j}(M_\infty, g_\infty).
\end{equation}
This contradicts the assumption that 
$\lim_{i \rightarrow \infty}  \lambda_{p,j}(M_i, g_i) \: = \: \infty$.

Now suppose that it is not true that there is a uniform lower bound
$a_{p,j}$ on  $\{\lambda_{p,j}(M, g)\}_{(M, g) \in {\cal M}}$
with the property that
$\lim_{j \rightarrow \infty} a_{p,j} \: = \: \infty$.
Then there are a number $C \: > \: 0$, a
sequence $\{ (M_i, g_i) \}_{i=1}^\infty$ in ${\cal M}$
and a sequence of integers $\{j_i\}_{i=1}^\infty$ such that
$\lim_{i \rightarrow \infty} j_i \: = \: \infty$ and for each $i$,
$\lambda_{p,j_i}(M_i, g_i) \: \le \: C$. Take a subsequence 
$\{ (M_i, g_i) \}_{i=1}^\infty$ and an
$(M_\infty, g_\infty)$ as before.
Then for each $j$,
\begin{equation} \label{3.5}
\lambda_{p,j}(M_\infty, g_\infty) \: \le \: \sup_{i \rightarrow \infty}
\lambda_{p,j_i}(M_\infty, g_\infty) \: \le \: 
\sup_{i \rightarrow \infty} e^{J \epsilon} \: 
\lambda_{p,j_i}(M_i, g_i) \: \le \: e^{J \epsilon} \:  C.
\end{equation}
This contradicts the discreteness of the spectrum of the $p$-form Laplacian
on $(M_\infty, g_\infty)$.
\end{pf}

Proposition \ref{Proposition 1} shows that in a certain sense, one has uniform
eigenvalue bounds in the noncollapsing case.  
It seems possible that for a given $n \in \Z^+$,
$K \in \R$ and $v, D \: > \: 0$, the collection ${\cal M}$ of connected
$n$-dimensional Riemannian manifolds $(M, g)$ with sectional curvatures
greater than $K$, volume greater than $v$ and diameter less than $D$
satisfies the hypotheses of Proposition \ref{Proposition 1}.  
It is known that there is
a finite number of homeomorphism types in ${\cal M}$
\cite{Grove-Petersen-Wu (1990)}. On the other hand, the analogous
space of metrics defined with Ricci curvature, instead of sectional curvature,
will generally not
satisfy the hypotheses of Proposition \ref{Proposition 1} 
\cite{Perelman (1997)}.

\end{document}